\documentclass[12pt]{amsart}
\usepackage{times}
\usepackage{amsfonts}
\usepackage{amsthm}
\usepackage{amsmath}
\advance \topmargin by -\headheight
\advance \topmargin by -\headsep
\evensidemargin \oddsidemargin
\newtheorem{theorem}{Theorem}[section]
\newtheorem{lemma}[theorem]{Lemma}
\newtheorem{corollary}[theorem]{Corollary}

\newtheorem{definition}[theorem]{Definition}

\begin{document}

\title{Graphs of functions and vanishing free entropy}

\author{Kenley Jung}

\email{kjung@math.ucla.edu}

\address{Department of Mathematics, University of California,
Los Angeles, CA 90095-1555,USA}

\begin{abstract} Suppose $X = \{x_1, \ldots, x_n\}$ is an $n$-tuple of selfadjoint elements in a tracial von Neumann algebra $M$.  If $z=z^* \in M$ and for some $y=y^*$ in the von Neumann algebra generated by $X$ $\delta_0(y, z) < \delta_0(y) + \delta_0(z)$, then $\chi(X \cup \{z\}) = -\infty$ (here $\chi$ and $\delta_0$ denote the microstates free entropy and free entropy dimension, respectively).   In particular, if $z$ lies in the von Neumann algebra generated by $X$, then $\chi(X \cup \{z\}) = -\infty$.  The statement and its proof are motivated by geometric-measure-theoretic results on graphs of functions.  A similar statement for the nonmicrostates free entropy is obtained under the much stronger hypothesis that $z$ lies in the algebra generated by $X$.  
\end{abstract} \maketitle

\section*{introduction}

A number of structure results in geometric measure theory have 
counterparts in Voiculescu's free probability.  In this correspondence the unit 
interval is replaced with an $n$-tuple generating a diffuse, injective von 
Neumann algebra, Cartesian products are replaced with freeness, the unit 
ball in $\mathbb R^n$ is replaced with a free family of $n$-semicirculars, 
Lebesgue/Hausdorff measure is replaced with free entropy/free Hausdorff 
entropy, and Minkowski dimension is replaced with free entropy dimension.  
The analogy provides new geometric-measure-theoretic tools for 
understanding free entropy.  For example, in \cite{jung:s1b}, 
Besicovitch's work on rectifiable $1$-sets was both the primary technical 
and heuristic motivation for introducing and studying the class of 
strongly $1$-bounded von Neumann algebras.  Another such example is the 
free entropy power inequality of \cite{dvv:entropypower} which can be 
recast in the form of the isoperimetric inequality involving 
$\epsilon$-neighborhoods 
and approximate perimeter; here the semicircular element replaces the unit 
ball and free additive convolution replaces the Minkowski sum (that equality is 
achieved exactly when the distribution is semicircular essentially follows 
from 
\cite{schultz}).

In this note I observe one more example of the connection 
between geometric measure theory and free probability.  Suppose that $X 
\subset \mathbb R^n$ is compact and $f: X \rightarrow \mathbb R$ is 
continuous.  Consider the graph $G(f)=\{(x, f(x)) : x \in X\} \subset 
\mathbb R^{n+1}$ of $f$. What can be said about the metric space structure 
of $G(f)$?  Clearly $G(f)$ is compact and its dimension (Hausdorff, 
Minkowski, or any fractal variant) is in the interval $[\dim X, \dim 
X+1]$.  While nothing more seems evident about $\dim G(f)$, at 
least $G(f)$ has Lebesgue measure zero (assuming $f \in L^{\infty}(X)$ 
will suffice).  This last fact about Lebesgue measure zero 
graphs has an analogue in free probability.

Assume $M$ is a tracial von Neumann algebra and $X = \{x_1,\ldots, x_n\}$ 
is a finite $n$-tuple of selfadjoint elements in $M$.  The 
$(m,k,\gamma)$-microstate space of $X$, $\Gamma(X;m,k,\gamma)$, is the set 
of all $n$-tuples of $k \times k$ selfadjoint matrices, $(a_1,\ldots, 
a_n)$, such that for any $1 \leq p \leq m$ and $1 \leq i_1, \ldots, i_p 
\leq n$

\[ |tr_k(a_{i_1} \cdots a_{i_p}) - \varphi(x_{i_1} \ldots x_{i_p})| < \gamma
\]

\noindent where $tr_k$ is the normalized trace on the $k \times k$ 
matrices.  The set of $n$-tuples of $k \times k$ selfadjoint matrices can 
be regarded as Euclidean space according to the norm-metric $|(a_1,\ldots, 
a_n)|_2 = (\sum_{i=1}^n tr_k(a_i^2))^{\frac{1}{2}}$.  Denote by '$\text{vol}$' Lebesgue 
measure with respect to this identification and successively define

\[ \chi(X;m,\gamma) = \limsup_{k \rightarrow \infty} \left [ k^{-2} \cdot \log (\text{vol}(\Gamma(X;m,k,\gamma)) + n \log k \right ],  \]

\[ \chi(X) = \inf \{ \chi(X;m,\gamma): m \in \mathbb N, \gamma >0 \}.
\]

\noindent  For the expert, notice that $n \log k$ appears instead of $\frac{n}{2} \log k$ due to the normalization of $\text{vol}$. $\chi(X)$ is called the free entropy of $X$ and can be regarded 
as the logarithmic Lebesgue measure of $X$.  There is also a Minkowski/box 
dimension quantity for $X$, $\delta_0(X)$, called the (modified) free 
entropy dimension of $X$.  For $\epsilon >0$, one replaces the volume of 
$\Gamma(X;m,k,\gamma)$ with the minimum number of $\epsilon$-balls 
required to cover $\Gamma(X;m,k,\gamma)$ in the above formulae and 
applying the same limits, arrives at the quantity $\mathbb 
K_{\epsilon}(X)$.  One then defines $\delta_0(X) = \limsup_{\epsilon 
\rightarrow 0} \frac{\mathbb K_{\epsilon}(X)}{|\log \epsilon|}$.

Suppose $y=y^*$ is in the von Neumann algebra generated by $X$.  Using the 
well-known heuristic that "$y= f(X)$ for $f \in L^{\infty}(X)$", it seems 
reasonable that the "graph of $f$", $X \cup \{y\} = X \cup \{f(X)\}$, 
should have "Lebesgue measure $0$", i.e., that its "logarithmic volume = 
free entropy" should be $-\infty$.  I will show that this is indeed the case.  
More generally, I will prove that if $z$ is a selfadjoint element such 
that $\delta_0(y, z) < \delta_0(y) + \delta_0(z)$, then $\chi(X \cup 
\{z\}) = -\infty$.  This is motivated by the fact that the graph of a 'typical'
continuous function from $\mathbb R$ into $\mathbb R$ will meet a non-orthogonal \emph{rotation} of the graph of a Lipschitz 
function in a set of positive Hausdorff $1$-measure (although its 
intersection with the graph of any Lipschitz function is a set of zero 
Hausdorff $1$-measure).  See \cite{mattila} and \cite{oneil:graph} for more 
details).

The proof bears some resemblance to the classical one.  There are, 
however, technical issues to resolve.  The tools for dispensing 
with these include the relative microstate decomposition in 
\cite{jung:hid} in both a free entropy and free entropy dimension form, 
estimates on intertwining sets found in \cite{gs}, covering number 
estimates of Borel subsets of the unitary group in terms of their actions 
on single selfadjoint matrices, and a covering estimate by relative 
microstate spaces.  For the rough outline of the proof and how it's 
related to the classical one, see the opening discussion of Section 2.

I will also prove a significantly weaker statement for the non-microstates 
free entropy $\chi^*$.  More specifically, I will show that if $y$ is in 
the \textit{algebra} generated by $X$, then $\delta^*(X \cup \{y\})   
<n+1$ and this will imply that $\chi^*(X \cup \{y\}) = -\infty$. 

There are 3 sections in this paper.  The first takes care of notation and 
reviews some facts which will be used later.  Section 2 proves the 
vanishing microstate entropy assertion stated in the abstract.  Also the 
equation, $\chi(F \cup \{a\}) = \chi(\Xi(F)) + \chi(a)$ 
and the inequality $\chi(F \cup G \cup \{a\}) \leq \chi(F \cup \{a\}) + 
\chi(G \cup \{a\}) - \chi(a)$ are established where $F$ and $G$ are tuples 
and $a$ is any selfadjoint element.  The last section concerns the 
nonmicrostate claim.

\section{preliminaries and notation}

In this section I'll set forth some notation and recall some facts that will be needed in the proofs.  

\vspace{.07in}

 Throughout this paper $M$ is a von Neumann algebra with a fixed tracial 
state $\varphi$.  $X =\{x_1,\ldots, x_n\}$ will denote an $n$-tuple of 
selfadjoint elements in $M$.  Given two $n$-tuples $Y= \{y_1,\ldots, 
y_n\}$ and $Z=\{z_1,\ldots, z_n\}$ in $M$, and $\epsilon >0$ define $Y + Z = 
\{y_1 + z_1,\ldots, y_n+ z_n\}$ and $\epsilon Y = \{ \epsilon y_1,\ldots, 
\epsilon y_n\}$.  Also, $W^*(Y)$ denotes the unital von Neumann algebra generated by $Y$.

\vspace{.07in}

For each $k$, $C_{k, n} = \frac{ \pi^{\frac{n 
k^2}{2}}}{\Gamma(\frac{nk^2}{2} +1)}$, the volume of the unit ball in 
$\mathbb R^{nk^2}$.  Stirling's Formula implies that $\lim_{k \rightarrow 
\infty} [ k^{-2} \cdot \log (C_{k,n})  - n \log k ] = \frac{n}{2}\cdot 
\log 
(\frac{2\pi}{n})$.

\vspace{.07in}
For each $k$ $U_k$ is the set of $k\times k$ unitaries, $M_k(\mathbb C)$ 
is the set of all $k\times k$ complex matrices, $M_k^{sa}(\mathbb C)$ is 
the set of all $k \times k$ selfadjoint complex matrices, and 
$(M^{sa}_k(\mathbb C))^n$ is the set of $n$-tuples of $k\times k$ 
selfadjoint matrices.  $tr_k$ is the normalized trace on the $k \times k$ 
complex matrices.  $|\cdot |_2$ is the Euclidean norm on 
$(M_k(\mathbb C))^n$ given by $|(a_1,\ldots, a_n)|^2_2 = \sum_{i=1}^n 
tr_k(a_j^* a_j)$.  Again, 'vol' denotes Lebesgue measure on $(M^{sa}_k(\mathbb C))^n$ with respect to this Euclidean identification via $| \cdot |_2$.  

\vspace{.1in}

For $G \subset (M_k(\mathbb C))^n$ $K_{\epsilon, \infty}(G)$ and $K_{\epsilon, 2}(G)$ denote the 
minimum number of open $\epsilon$-balls of $G$ required to cover $G$ with respect 
to the operator norm and $|\cdot 
|_2$-norm, respectively.   Similarly $P_{\epsilon, \infty}(G)$ and $P_{\epsilon, 2}(G)$ denote the maximum number of elements in a collection of mutually disjoint open $\epsilon$-balls of $G$ with respect to the operator and $|\cdot |_2$-norm, respectively.  As is well-known (\cite{dvv:entropysurvey}), one can 
interchange these quanities with impunity when it comes to their 
asymptotic behavior.

\vspace{.07in}

Recall the covering estimates for intertwining sets in \cite{gs}.  There 
the authors considered for each $k$, $z_k$, the diagonal unitary whose 
entries list in counterclockwise order, the $k$ roots of unity.  They then 
considered the set $\Sigma(z_k, z_k, \epsilon) = \{ a \in M_k(\mathbb C): \|a\| \leq 2,
|a z_k - z_k a|_2 < \epsilon\}$ and showed that 

\[ K_{\epsilon, 2}(\Sigma(z_k,z_k, \epsilon^2)) \leq \left ( 
\frac{36}{\epsilon} \right)^{4 \pi k^2 \epsilon + 2k}. \]

\noindent They actually showed this in greater generality but this is not 
important for the purposes at hand.  The only important constants are the 
exponent on the left hand side.

\vspace{.1in}

Suppose $y$ is a selfadjoint element in 
$M$ such that the distribution of $y$ induced by $\varphi$ is Lebesgue 
measure on the unit interval.  Observe that $\chi(y) > -\infty$.  For each 
$k$ denote by $y_k$ the $k \times k$ diagonal matrix whose entries are $0, 
\frac{1}{k}, \frac{2}{k}, \ldots, \frac{k-1}{k}$, listed along the 
diagonal in increasing order.  Notice that for any $m \in \mathbb N$ and 
$\gamma >0$ $y_k \in \Gamma(y;m,k,\gamma)$ for $k$ sufficiently large.  
Define $\mathcal I(y_k, y_k, \epsilon) = \{ u \in U_k: |u y_k - y_k u|_2 < 
\epsilon \}$.  If $u \in \mathcal I(y_k, y_k, \epsilon)$, then $|u y_k u^* 
- y_k|_2 < \epsilon$ which implies

\begin{eqnarray*} |u z_k - z_k u|_2 & = & |u z_k u^* - z_k|_2 \\ & = & |u 
e^{2\pi i y_k} u^* - e^{2 \pi i y_k}|_2 \\ & \leq & 2 \pi |u y_k 
u^* - y_k|_2 \\ & \leq & 2 \pi \epsilon. 
\end{eqnarray*}

\noindent From this it follows that $\mathcal I(y_k, y_k, \epsilon^2) 
\subset \Sigma(z_k, z_k, 7 \epsilon^2)$ so that,

\[ K_{\epsilon, 2}(\mathcal I(y_k, y_k, \epsilon^2)) \leq 
\left( \frac{252}{\epsilon} \right)^{28 \pi k^2 
\epsilon + 2k}.\]

\noindent It is a well-known fact that the volume of 
the operator norm ball of $M^{sa}_k(\mathbb C)$ of 
radius $1$ and the $|\cdot|_2$ norm ball of 
$M^{sa}_k(\mathbb C)$ of radius $1$ are exponentially 
proportional in $k^2$ (\cite{raymond}), and thus there 
exists a constant $C >0$ independent of $k$ such that 
for $\epsilon \in (0,1/2)$

\[ K_{\epsilon, \infty}(\mathcal I(y_k, y_k, \epsilon^2)) 
\leq \left( \frac{1}{\epsilon}\right)^{Ck^2 
\epsilon}.\]

\vspace{.07in}

Finally, recall from \cite{szarek} and \cite{vn} that there exist constants $\alpha_0, \alpha_1  >0$ such that for any $\epsilon \in (0,1)$,

\[ \left (\frac{\alpha_0}{\epsilon} \right)^{k^2} \leq K_{\epsilon, \infty}(U_k) \leq  \left(\frac{\alpha_1}{\epsilon} \right)^{k^2}.\]

\section{vanishing free entropy}

Before I begin the proof consider how one would show that the graph 
of a bounded, measurable real-valued function $f$ from $[0,1]$ has $0$ 
area.  By scaling and translating, assume without loss of generality that 
$f$ takes values in $[0,1]$.  For $\epsilon >0$ partition $[0,1]$ into $n$ 
mutually disjoint sets $I_1, \ldots, I_n$ each of length no greater than 
$\epsilon$.  

\begin{eqnarray} G(f) \subset \cup_{j=1}^n f^{-1}(I_j) \times I_j
\end{eqnarray}

\noindent so that

\[ \text{vol}(G(f)) \leq  \sum_{j=1}^n \epsilon \cdot f^{-1}(I_j) =
\epsilon.
\]

Turning to the free entropy situation, consider the problem of showing $\chi(X \cup \{y\}) = 
-\infty$ where $y=y^*$ is an element in the von Neumann algebra 
generated by $X$ with uniform distribution on $[0,1]$ (this is not as general as the stated result in the introduction, but a good portion of the issues are translucent in this easier setting).  To estimate 
$\chi(X \cup \{y\})$ one needs a handle on the microstate space $\Gamma(X 
\cup y;m,k,\gamma)$. An arbitrary microstate space has an obvious symmetry: 
it is invariant under unitary conjugation.  This allows one to write the 
space as a product of $U_k$ with an appropriate quotient space.  
Formally, for any $m, k \in \mathbb N$ and $\gamma
>0$ define

\[ \Xi(X;m,k,\gamma) = \{ \xi \in (M^{sa}_k(\mathbb C))^n: (\xi, y_k) \in 
\Gamma(X \cup 
\{y\}:m,k,\gamma)\}.
\]

\noindent Resorting to the heuristic 
"$f(X) = y$", write $\Xi(X;m,k,\gamma) = f^{-1}(y_k)$.  Now, any 
microstate of $y$ can be 
approximated by a unitary conjugate of $y_k$ (since $y$ 
generates an abelian, and thus injective von Neumann algebra) and because 
$\Gamma(X \cup 
\{y\};m,k,\gamma)$ is invariant under the $U_k$-action, it follows that

\begin{eqnarray} \Gamma(X \cup \{y\};m,k,\gamma) \subset \cup_{u \in U_k} 
u (f^{-1}(y_k) \times \{y_k\}) u^*.
\end{eqnarray}

\noindent Notice the similarity between (1) and (2) and that in (2) the 
sets in the union are unitarily conjugate, therefore isometric, and thus 
have the same volume.  Using this fact it is now tempting to "sum" or 
integrate the right hand side as in the classical case.  For the 
microstate setting it will be more convenient to use packings by balls and 
by isometric copies of $f^{-1}(y_k)$.  Using packings, I will argue that the unitary action on the $f^{-1}(y_k) \times \{y_k\}$ is "approximately faithful" and thus, the union in the set above is almost disjoint, so that for any $\epsilon >0$, the volume of the union is almost the product of the $\epsilon$-metric entropy of $U_k$ times the volume of $f^{-1}(y_k) \times \{y_k\}$.  Since this union is contained in a bounded ball in $(M^{sa}_k(\mathbb C))^n$, this will give an appropriate bound on the volume of $f^{-1}(y_k) \times \{y_k\}$, and ultimately on $\chi(X \cup \{y\})$.  

Now for the technicalities (=proof).  The first thing I'll do is prove a 
relative microstates entropy equation.  This is essentially an entropy 
version of the relative free entropy dimension equation in 
\cite{jung:hid}.  Suppose $a \in M$ is a selfadjoint element in $M$ and 
fix a sequence $\langle a_k \rangle_{k=1}^{\infty}$ such that for any $m 
\in \mathbb N$ and $\gamma
>0$, $a_k \in \Gamma(a;m,k,\gamma)$ for $k$ sufficiently large and such that the operator norms of the $a_k$ are uniformly bounded.  Assume 
$F$ is an $n$-tuple of selfadjoint elements in $M$. Denote by $\Xi(F; 
\cdot)$ all the microstate quantities takes relative to this sequence 
$\langle a_k \rangle_{k=1}^{\infty}$, i.e., for any $m, k \in \mathbb N$ 
and $\gamma
>0$

\[ \Xi(F;m,k,\gamma) = \{ \xi \in (M^{sa}_k(\mathbb C))^n: (\xi, a_k) \in 
\Gamma(F \cup 
\{a\}:m,k,\gamma)\}. \]

\noindent The use of an operator norm restriction will be convenient.  
For any $R >0$ define $\Xi_R(F;m,k,\gamma)$ to be the set of 
all $n$-tuples $\xi$ such that $\xi \in 
\Gamma(F;m,k,\gamma)$ and such that every entry in $\xi$ has operator norm 
no greater than $R$.

\begin{definition} Given $a \in M$ and $\langle a_k 
\rangle^{\infty}_{k=1}$ 
as above, define successively,

\[ \chi(\Xi(F;m,\gamma)) = \limsup_{k \rightarrow \infty} \left [k^{-2} 
\cdot \log (\text{vol}(\Xi(F;m,k,\gamma)) + n \log k \right], 
\]

\[ \chi(\Xi(F)) = \inf \{ \chi(\Xi(F;m,\gamma)): m \in \mathbb N, \gamma >0\}.
\]

\noindent Also for any $R >0$, define 
$\chi_R(\Xi(F;m,\gamma))$ and $\chi_R(\Xi(F))$ to be 
the resultant quantities obtained by replacing 
$\Xi(F;m,k,\gamma)$ in the definition above with 
$\Xi_R(F;m,k,\gamma)$. \end{definition}

The above definition of $\chi(\Xi(F))$ may seem 
dependent on $a_k$.  It doesn't.  First observe that if 
$R$ is strictly greater than the operator norms of any 
of the elements in $F \cup \{a\}$, then $\chi(\Xi(F)) = 
\chi_R(\Xi(F))$.  Clearly $\chi(\Xi(F)) \geq 
\chi_R(\Xi(F))$.  For the reverse inequality the lower 
bounds in \cite{sh} (see also \cite{dvv:entropy2}) for 
the Jacobian of a coordinatewise spectral cutoff 
function from $\Gamma(F \cup \{a\};m,k,\gamma)$ into 
$\Gamma_R(F \cup \{a\};m^{\prime},k,\gamma^{\prime})$ 
imply that $\chi(\Xi(F)) \leq \chi_R(\Xi(F))$.  Thus, 
$\chi(\Xi(F)) = \chi_R(\Xi(F))$.

Now suppose $\langle b_k \rangle_{k=1}^{\infty}$ is 
another sequence such that for any $m \in \mathbb N$ 
and $\gamma >0$, $b_k \in \Gamma(y;m,k,\gamma)$ for $k$ 
sufficiently large and such that the operator norms of 
the $b_k$ are uniformly bounded.  Denote the relative 
microstate quantities of $F$ taken with respect to the 
sequence $\langle b_k \rangle_{k=1}^{\infty}$ by 
$\Xi^{\prime}(F;m,k,\gamma)$, 
$\Xi_R^{\prime}(F;m,k,\gamma)$ $\chi(\Xi^{\prime}(F))$, 
and $\chi(\Xi_R^{\prime}(F))$.  I claim that 
$\chi(\Xi(F)) = \chi(\Xi^{\prime}(F))$.  In order to 
prove this it suffices by the preceding paragraph to 
show that $\chi_R(\Xi(F)) = \chi_R(\Xi^{\prime}(F))$ 
for some $R$ greater than the operator norms in $F \cup 
\{a\}$.  Choose such an $R$ and make sure it is greater 
than the operator norms of the $a_k$ and $b_k$.  
Suppose $m_0 \in \mathbb N$ and $\gamma_0 >0$. There 
exist a $m_0 \leq m \in \mathbb N$ and $\gamma_0 > t, \gamma >0$ such that 
if $\xi \in \Gamma_R(F \cup \{a\};m,k,\gamma)$ and 
$|\eta - \xi|_2 <t$ and the entries of $\eta$ have 
operator norms no greater than $R$, then $\eta \in 
\Gamma_R(F \cup \{a\};m_0,k,\gamma_0)$.  There also 
exists a sequence $\langle u_k \rangle_{k=1}^{\infty}$ 
in $U_k$ such that $|u_k a_k u_k^* - b_k |_2 
\rightarrow 0$ and from this it follows that for 
sufficiently large $k$

\[ u_k \Xi_R(F;m,k,\gamma)u_k^* \subset \Xi_R^{\prime}(F;m_0,k,\gamma_0).
\]

\noindent As $m_0$ and $\gamma_0$ were arbitrary and 
unitary conjugation is an isometric, and thus measure 
preserving action, it follows that $\chi(\Xi_R(F)) \leq 
\chi_R(\Xi^{\prime}(F))$.  The reverse inequality 
follows from the exact same argument and so 
$\chi_R(\Xi(F)) = \chi_R(\Xi^{\prime}(F))$.  By 
taking 
$R$ sufficiently large it follows that $\chi(\Xi(F)) = 
\chi^{\prime}(\Xi(F))$.

 The following now makes sense:

\begin{definition} Suppose $\langle a_k 
\rangle_{k=1}^{\infty}$ is a sequence such that for any 
$m \in \mathbb N$ and $\gamma >0$, $a_k \in 
\Gamma(a;m,k,\gamma)$ for $k$ sufficiently large and 
such that the sequence is uniformly bounded in operator 
norm.  The free entropy of $F$ relative to $a$, denoted 
by $\chi(F/a)$, is $\chi(\Xi(F))$ where this latter 
quantity is computed according to Definition 2.1 and 
with respect to the sequence $\langle a_k 
\rangle_{k=1}^{\infty}$.  If $R >0$ is strictly greater 
than the operator norms of any of the elements in $F$ 
and $\langle a_k \rangle_{k=1}^{\infty}$ then 
$\chi(F/a) = \chi_R(\Xi(F))$ where this latter 
quantity 
again, is computed according to Definition 2.1 with 
respect to the sequence $\langle a_k 
\rangle_{k=1}^{\infty}$.  \end{definition}

Now for the promised relative entropy formula:

\begin{lemma} $\chi(F \cup \{a\}) = \chi(F/a) + \chi(a)$.
\end{lemma}

\begin{proof} Find and fix a sequence $\langle a_k 
\rangle_{k=1}^{\infty}$ such that for any $m \in 
\mathbb N$ and $\gamma >0$, $a_k \in 
\Gamma(a;m,k,\gamma)$ for $k$ sufficiently large and 
such that the sequence is uniformly bounded in operator 
norm and consider all relative microstate spaces with 
respect to this sequence $\langle a_k 
\rangle_{k=1}^{\infty}$.  Find an $R$ strictly greater 
than the operator norm of any element in $F \cup \{a\} 
\cup \langle a_k \rangle_{k=1}^{\infty}$.  By \cite{sh} and Definition 2.2, it suffices 
to show $\chi_R(F \cup \{a\}) = \chi_R(\Xi(F)) + 
\chi(a)$.

Let $m \in \mathbb N$ and $\gamma >0$ be given.  There exist $t, \gamma_0 >0$ so small that if $(\xi, x) \in \Gamma_R(F \cup \{a\}; m,k, \gamma_0)$, $|x- y|_2 <t$, and $\|y\| \leq R$, then $(\xi, y) \in \Gamma(F \cup \{a\};m,k,\gamma)$.  Find $m_1 \in \mathbb N$ and $\gamma_1>0$ such that if $x, y \in \Gamma(a;m_1, k,\gamma_1)$, then there exists a $u \in U_k$ such that $|uxu^* -y|_2 < t/2$.   Now suppose $m_2 =m +m_1$ and $\gamma_2 = \min \{ \gamma_0, \gamma_1\}$.  For each $x \in \Gamma_R(a;m_2,k,\gamma_2)$ define $L_x = \{\xi \in (M^{sa}_k(\mathbb C))^n : (\xi, x) \in \Gamma(X \cup \{a\};m_2, k, \gamma_2)\}$ and $C_x = \{\xi \in (M^{sa}_k(\mathbb C))^n : (\xi, x) \in \Gamma_R(X \cup \{a\};m,k,\gamma)\}$.  By Fubini's Theorem

\begin{eqnarray}
\int_{\Gamma_R(a;m_2,k,\gamma_2)} \text{vol}(L_x) \, dx = \text{vol}\left(\Gamma_R(F \cup \{a\};m_2,k,\gamma_2 )\right).
\end{eqnarray}

\noindent and

\begin{eqnarray} \text{vol}\left(\Gamma_R(F \cup \{a\};m,k,\gamma) \right) \geq  \int_{\Gamma_R(a;m_2,k,\gamma_2)} \text{vol}(C_x) \, dx 
\end{eqnarray}

\noindent where $dx$ is Lebesgue measure on $M^{sa}_k(\mathbb C)$.  For $k$ sufficiently large and each $x \in \Gamma_R(a;m_2,k,\gamma_2)$ there exists a $u \in U_k$ such that $|uxu^* - a_k|_2 < t/2$ which implies that for any $\xi \in L_x$, $(u \xi u^*, a_k ) \in \Gamma_R(F \cup \{a\};m,k,\gamma)$, i.e., $uL_x u^* \subset \Xi(F;m,k,\gamma) \Rightarrow \text{vol}(L_x) = \text{vol}(uL_xu^*) \leq \text{vol}(\Xi_R(F;m,k,\gamma))$.  Thus, (3) implies 

\[ \text{vol}\left(\Gamma_R(F \cup \{a\};m_2, k,\gamma_2) \right) \leq  \text{vol}(\Xi_R(F;m,k,\gamma)) \cdot   \text{vol}\left(\Gamma_R(a;m_2,k,\gamma_2)\right).
\]

\noindent Similarly $|x - u^*a_ku|_2 < t/2$ implies that for any $\xi \in L_{ua_ku^*}$, $(\xi, x) \in \Gamma_R(F \cup\{a\};m,k,\gamma)$, i.e., $L_{ua_ku^*} \subset C_x \Rightarrow \text{vol}(\Xi(F;m_2,k,\gamma_2)) = \text{vol}(L_{ua_ku^*}) \leq \text{vol}(C_x)$.  Thus, (4) implies 

\[ \text{vol}\left(\Gamma_R(F \cup \{y\};m,k,\gamma) \right) \geq \text{vol}(\Xi_R(F;m_2,k,\gamma_2)) \cdot \text{vol}\left(\Gamma_R(a;m_2,k,\gamma_2) \right).
\]

  For any $m \in \mathbb N$ and $\gamma >0$ there exist $m_2 >m$ and $0 < \gamma_2 < \gamma$ such that 

\[ \chi_R(F \cup \{a\};m_2,\gamma_2) \leq 
\chi_R(\Xi(F;m,\gamma)) + \chi_R(a;m_2, \gamma_2)
\]

\noindent and

\[ \chi_R(\Xi(F;m_2,\gamma_2) + \chi_R(a;m_2,\gamma_2) 
\leq \chi_R(F \cup \{a\};m,\gamma).
\]

\noindent From this it follows that $\chi_R(F \cup 
\{a\}) = \chi_R(\Xi(F)) + \chi_R(a)$. \end{proof}

Before going on to the main argument observe the following which is the free entropy version of the hyperfinite inequality for $\delta_0$ in \cite{jung:hid}

\begin{corollary} If $F, G$ are finite tuples in $M$ and $a=a^* \in M$ 
such that $\chi(a) > - \infty$, 
then 

\[\chi(F \cup G \cup \{a\}) \leq \chi(F \cup \{a\}) + \chi(G \cup 
\{a\}) - \chi(a).\]
\end{corollary}

\begin{proof} $-\infty< \chi(a) < \infty$ and $\chi((F \cup G)/a) \leq 
\chi(F/a) + \chi(G/a)$ since $\Xi(F \cup G;m,k,\gamma) \subset 
\Xi(F;m,k,\gamma) \times \Xi(G;m,k,\gamma)$.  Using Lemma 2.3 and keeping 
in mind that all quantities below lie in $[-\infty, \infty)$, 

\begin{eqnarray*} \chi(F \cup G \cup \{a\}) & = & \chi((F \cup G)/a) + \chi(a) \\ & \leq & \chi(F/a) + \chi(G/a) + \chi(a) \\ & = & \chi(F/a) + \chi(G/a) + 2 \chi(a) - \chi(a) \\ & = & \chi(F) + \chi(G) - \chi(a).
\end{eqnarray*}
\end{proof}

In the remainder of this section $y_k$ is the $k \times k$ diagonal matrix whose $ii$th entry is $\frac{i-1}{k}$.  Notice that if $y = y^* \in M$ is any element whose distribution with respect to $\varphi$ is Lebesgue measure on the $[0,1]$, then for any $m \in \mathbb N$ and $\gamma >0$, $y_k \in \Gamma(y;m,k,\gamma)$ for $k$ sufficiently large.  The next lemma I'll need provides an upper bound for the asymptotic covering numbers of Borel subsets of $U_k$, provided such bounds exist for the spaces obtained from the action of the Borel sets on the $y_k$.  

Recall the following well-known random matrix result.  Denote by $H_k$ the 
group of diagonal unitaries and by $\mathbb R^k_<$ the set of all 
$(t_1,\ldots, t_k) \in \mathbb R^k$ such that $t_1 < \cdots < t_k$.  
There exists a map $\Phi_k: M^{sa}_k(\mathbb C) \rightarrow U_k/H_k \times 
\mathbb R^k_<$ defined almost everywhere on $M^{sa}_k(\mathbb C)$ such 
that for each $x \in M^{sa}_k(\mathbb C)$, $\Phi_k(x) = (h,z)$ where $z$ 
lists the eigenvalues of $x$ in increasing order and $h$ is the image of 
any unitary $u$ in $U_k/G_k$ satisfying $uzu^* =x$.  Notice that $\Phi_k$ 
is a bijection from the set of all selfadjoint $k\times k$ matrices 
with $k$ 
distinct eigenvalues onto $U_k/H_k \times \mathbb R^k_<$.  The pushforward 
of vol by the map $\Phi_k$ induces 
a measure $\mu_k$ on $U_k/H_k \times \mathbb R^k_<$ and moreover,

\[\mu_k = \nu_k \times D_k \int_{\mathbb R^k_<} \Pi_{1 \leq i < j \leq k} (t_i -t_j) \, dt_1 \cdots dt_k,\]

\noindent where $D_k = \frac{(2\pi)^{\frac{k(k-1)}{2}}}{\sqrt{k} 
\Pi_{j=1}^k j!}$ and $\nu_k$ is the probability measure on $U_k/H_k$ 
induced by Haar probability measure.  See \cite{mehta} for more details.

\begin{lemma}  Suppose for each $k$, $G_k$ is a Borel subset of $U_k$ such that 
for 
some fixed $\beta \in (0,1)$ and all $\epsilon >0$ sufficiently small 

\[ \beta |\log \epsilon| \geq \limsup_{k \rightarrow \infty} k^{-2} \cdot \log 
\left [K_{\epsilon, 2}(\{uy_ku^*: u \in G_k\}) \right].
\]

\noindent Then there exists an $\epsilon_0 >0 $ such that for $0 < 
\epsilon <\epsilon_0$ 

\[\beta |\log \epsilon| - \chi(y) + \log \alpha_1 + \log 5 > 
\limsup_{k 
\rightarrow \infty} k^{-2} \cdot \log (K_{\epsilon, \infty}( G_k)).
\]
\end{lemma}

\begin{proof} The $3 \epsilon$-neighborhood of $\{uy_k u^*: u \in G_k\} \subset M^{sa}_k(\mathbb C)$
certainly contains all elements of the form $vav^*$ where $v \in U_k$ and 
$\| v-u\| < \epsilon$ for some $u \in G_k$ and $a$ is a $k \times k$ 
matrix in $\Omega_k$ where $\Omega_k$ consists of all diagonal matrices whose eigenvalues listed in increasing order lie in the set
 
\[  R_k =[- 1/4k, 1/4k ] \times [1/k - 1/4k, 1/k 
+ 1/4k] \times \cdots \times [(k-1)/k - 1/4k, (k-1)/k + 1/4k]. \]

\noindent Denoting by $\pi_k: U_k \rightarrow U_k/H_k$ the quotient map, the discussion preceding the lemma implies

\begin{eqnarray*} K_{\epsilon, 2}(\{u y_k u^*: u \in G_k\}) \cdot C_k 
\cdot (4\epsilon)^{k^2} & \geq &\text{vol}\left( \mathcal 
N_{3\epsilon}(\{uy_k u^*: u \in G_k\}) \right) \\ & \geq & 
\text{vol}(\{uyu^*: u \in \mathcal N_{\epsilon, \infty}(G_k), y \in 
\Omega_k\}) \\ & = & \mu_k(\Phi(\{uyu^*: u \in \mathcal N_{\epsilon, 
\infty}(G_k), y\in \Omega_k\}) \\ & = & \nu_k(\pi_k(\mathcal N_{\epsilon, 
\infty}(G_k))) \cdot D_k \cdot \int_{R_k} \Pi_{1 \leq i < j \leq k} (t_i - 
t_j)^2 \, dt_1 \cdots dt_k.\\ \end{eqnarray*}

\noindent Computing exactly as in Proposition 4.5 of \cite{dvv:entropy2}

\[\int_{R_k} \Pi_{1 \leq i < j \leq k} (t_i -t_j)^2 \, dt_1 \cdots dt_k \geq  (4k)^{-1} \cdot \Pi_{1 \leq i < j \leq k} \left ( \frac{(j-i)}{k^2} - \frac{1}{2k}\right)^2
\]

\noindent so that $\liminf_{k \rightarrow \infty} k^{-2} \cdot \log \left[ C_k^{-1} \cdot D_k \cdot \int_{\Omega_k} \Pi_{1 \leq i < j \leq k} (t_i - t_j)^2 \, 
dt_1 \cdots dt_k \right ]$ dominates

\begin{eqnarray*} 2^{-1} \log 2\pi + \frac{3}{4} + \liminf_{k \rightarrow \infty} k^{-2} \cdot \sum_{1 \leq i \leq j \leq k} \log \left (\frac{(j -i)}{k} - \frac{1}{2k} \right )^2= \chi(y) > -\infty.\\
\end{eqnarray*}

\noindent Thus,

\[ \limsup_{k \rightarrow \infty} k^{-2} \cdot \log \left[ K_{\epsilon, 2}(\{uy_ku^*: u \in G_k\} ) (4 \epsilon)^{k^2} \right] \geq \limsup_{k \rightarrow \infty}[\nu_k(\pi_k(\mathcal N_{\epsilon, \infty}(G_k)))] + \chi(y).
\] 

Combining this last estimate above with the hypothesis of the lemma,

\begin{eqnarray*} (\beta - 1) |\log \epsilon| & \geq & \limsup_{k 
\rightarrow \infty} k^{-2} \cdot \log \left [ K_{\epsilon,2}(\{uy_k u^*: u 
\in G_k\}) \cdot (4 \epsilon)^{k^2} \right ] - \log 4 \\ & \geq & 
\limsup_{k \rightarrow \infty} \left [k^{-2} \cdot \log \nu_k (\pi_k 
(\mathcal N_{\epsilon, \infty}(G_k))) \right] + \chi(y) - \log 4. \\ 
\end{eqnarray*}

\noindent  Denote by $m_k$ normalized Haar measure on $U_k$ and write $B_{\epsilon}$ for the open ball with center the identity operator and radius $\epsilon$ with respect to the operator norm.   By definition of $\nu_k$,

\begin{eqnarray*} \nu_k(\pi_k(\mathcal N_{\epsilon, \infty}(G_k))) & \geq & m_k(\mathcal N_{\epsilon, \infty}(G_k)) \\ & \geq & P_{\epsilon, \infty}(G_k) \cdot m_k(B_{\epsilon}) \\ & \geq & K_{2\epsilon, \infty}(G_k) \cdot m_k(B_{\epsilon}) \\ & \geq & K_{2 \epsilon, \infty}(G_k) \cdot P_{\epsilon, \infty}(U_k)^{-1} \\ & \geq & K_{2 \epsilon, \infty}(G_k) \cdot \epsilon^{k^2} \alpha_1^{-k^2}.
\end{eqnarray*}

\noindent Substituting this into the previous inequalities,

\[ (\beta - 1)|\log \epsilon| \geq \limsup_{k \rightarrow \infty} \left [ 
k^{-2} \cdot \log (K_{2 \epsilon, \infty}(G_k)) \right ] + \log \epsilon - 
\log \alpha_1 + \chi(y) - \log 4 \]

\noindent which gives the desired claim.
\end{proof}

\begin{theorem} If $x=x^*$ is an element in the von Neumann algebra 
generated by $X$ and $z =z^* \in M$ such that $\delta_0(x, z) < 
\delta_0(x) + \delta_0(z)$, then $\chi(X \cup \{z\}) = - \infty$.
\end{theorem}

\begin{proof} I'll break this proof up into several parts.

\vspace{.1in}

\noindent Part 1: First I'll reduce the problem to the case where there exists 
a $y=y^*$ in the von Neumann algebra generated by $X$ such that 
$\delta_0(y,z) <2$ and the distribution of $y$ induced by $\varphi$ is 
Lebesgue measure on the unit interval.  Without loss of generality 
$X$ generates a diffuse von Neumann algebra and $\delta_0(z) = 1$.  Find a 
$y = y^*$ in the von Neumann algebra generated by $X$ such that the 
distribution of $y$ induced by $\varphi$ is Lebesgue measure on the unit 
interval, and such that $x$ is contained in the von Neumann algebra 
generated by $y$. By \cite{jung:hid}

\begin{eqnarray*} \delta_0(y, z) & = & \delta_0(x,y,z) \\ & \leq & 
\delta_0(x, z) + \delta_0(x,y) - \delta_0(x) \\ & < & \delta_0(x) + 
\delta_0(z) + \delta_0(y) - \delta_0(x) \\ & = & 2. \end{eqnarray*}

\vspace{.1in}

\noindent Part 2: By Part 1 there exists 
a $y = y^*$ in the von Neumann algebra generated by $X$ such that 
$\delta_0(y, z) <2$ and the distribution of $y$ induced by $\varphi$ is 
Lebesgue measure on the unit interval.  I want to show that $\chi(X \cup 
\{z\}) = - \infty$, and by Lemma 2.3 it suffices to show that $\chi(X/z) = 
-\infty$. Find an $R >1$ greater than the operator norms of any of the 
elements in $X \cup \{z\}$.  Pick a sequence $\langle z_k 
\rangle_{k=1}^{\infty}$ such that for any $m \in \mathbb N$ and $\gamma
>0$, $z_k \in \Gamma_R(z;m,k,\gamma)$ for $k$ sufficiently large.  
Consider all relative microstate spaces $\Xi()$ with respect to this 
sequence $\langle z_k \rangle_{k=1}^{\infty}$.  What I will do in this part is bound $\chi_R(\Xi(X))$ in terms of sets obtained from Borel subsets of $U_k$ and the relative microstates of $X$ with respect to a well approximating microstate sequence for $y$.  The Borel subsets will have an upper bound on their metric entropy in terms of $\delta_0(y, z) -1 <1$.    

The relative microstates decomposition for $\delta_0$ in \cite{jung:hid} yields

\[ 1 > \delta_0 (y, z) - \delta_0(z) = \limsup_{\epsilon \rightarrow 0} 
\frac{\mathbb K_{\epsilon}(\Xi(y))}{|\log \epsilon|}.
\]

\noindent There exist $1 > \beta, \epsilon_0 >0$ such that for all 
$\min\{\epsilon_0, 1/2\} > \epsilon >0$,

\[ \beta \cdot |\log \epsilon| \geq \mathbb K_{\epsilon}(\Xi(y)).
\]

Suppose 
$\epsilon_0 > \epsilon >0$.  Choose $m_1, \gamma_1$ (dependent on 
$\epsilon$) such that

\begin{eqnarray} \beta \cdot |\log \epsilon| \geq \limsup_{k \rightarrow 
\infty} k^{-2} 
\cdot \log (K_{\epsilon,2}(\Xi_R(y;m_1,k,\gamma_1)))
\end{eqnarray}

\noindent and such that there exists a selfadjoint polynomial $f$ in $n$ 
noncommuting variables satisfying the condition that for any $(\xi, a) \in 
\Gamma(X \cup \{y\};m_1,k,\gamma_1)$, $|f(\xi)
- a|_2 < \epsilon^2/2$.  There exist $t, \gamma_0 >0$ such that if 
  $(\xi, a) \in
\Xi_R(X \cup \{z\};m_1, k,\gamma_0)$, $\|b\| \leq R$, and $|a-b|_2<t$, 
then $(\xi, b) \in \Xi_R(X \cup \{y\};m_1,k,\gamma_1)$.  There also exists 
$m_2$, $\gamma_2$ such that $m_2 > m_1$, $\gamma_0 > \gamma_2$, and if $a, b \in \Gamma_R(z;m_2, k, \gamma_2)$, 
then there exists a $u \in U_k$ such that $|uau^* - b|_2 <t$.  Find a 
polynomial $g$ in $n$ noncommuting variables such that for some $m_3 \in 
\mathbb N$, $\gamma_3 >0$ and any $\xi \in \Xi_R(X;m_3,k,\gamma_3)$ (see 
Lemma 4.1 of \cite{jung:hf}),

\[ (\xi, g(\xi)) \in \Xi_R(X \cup \{y\};m_2, k, \gamma_2) \subset \Xi_R(X \cup\{y\};m_1,k,\gamma_0). 
\]

\noindent For $k$ sufficiently large there exists a $u \in U_k$ such 
that $|g(\xi) - u y_k u^*|_2 <t$ and hence,

\[ (\xi, uy_k u^*) \in \Xi_R(X \cup y;m_1,k,\gamma_1) \subset 
\Xi_R(X;m_1k,\gamma_1) \times \Xi_R(y;m_1,k,\gamma_1).
\]

\noindent The above is true for any $\xi \in \Xi_R(X;m_3,k,\gamma_3)$. 
Denote by $G_k$ the set of all $u \in U_k$ such that $u y_k u^* \in 
\Xi_R(y;m_1, k,\gamma_1)$ and set

\[ L(m_1,k,\gamma_1) = \{\xi \in (M^{sa}_k(\mathbb 
C))^n: (\xi, y_k) \in 
\Gamma_R(X \cup\{y\};m_1,k,\gamma_1)\}.\] 

\noindent $L(m_1,k,\gamma_1)$ is just the set of 
relative microstates of $X$ with respect to the 
sequence 
$\langle y_k \rangle_{k=1}^{\infty}$.  From the above
 
\[\Xi_R(X;m_3,k,\gamma_3) \subset \cup_{u \in G_k} u L(m_1,k,\gamma_1) 
u^*. \]

\noindent Set $F_k = \cup_{u \in G_k} u L(m_1,k,\gamma_1) u^*$.  

\vspace{.1in}

\noindent Part 3: I now want find an upper bound for $\text{vol} (F_k)$; 
by Part 2 this will serve as an upper bound for 
$\text{vol}(\Xi_R(X;m_3,k,\gamma_3))$. It will be easier to first find an 
upper bound for the volume of $\Omega_k = \cup_{u \in B_k} 
uL(m_1,k,\gamma_1) u^*$ where $B_k \subset U_k$ is the unit ball of 
operator norm radius $\epsilon$, centered at the identity.  That is what I 
will do in this part of the proof.

If 
$u L(m_1,k,\gamma_1) u^* \cap 
v L(m_1,k,\gamma_1) v^* \neq \emptyset$, then there exists a $\xi
\in L(m_1,k,\gamma_1)$ with $\xi= 
u^* v \eta v^* u$ and $\eta \in L(m_1,k,\gamma_1)$.  

\begin{eqnarray*} |u^* v y_k v^* u - y_k|_2 & \leq & |u^* v y_k v^* u - 
u^* v f(\eta) v^* u|_2 + |f(\xi) - y_k|_2 \\ & < & \epsilon^2. 
\end{eqnarray*}

\noindent Hence $u^* v \in \mathcal I(y_k, y_k,  
\epsilon^2) \iff u \in v [\mathcal I (y_k, y_k, \epsilon^2)]^*$.  So

\begin{eqnarray} u L(m_1,k,\gamma_1)u^* \cap v L(m_1,k,\gamma_1) v^* \neq 
\emptyset \Rightarrow u \in v[\mathcal I (y_k, y_k, \epsilon^2)]^*. 
\end{eqnarray}

 By the discussion in Section 1, 
fix an $\epsilon$-cover of unitaries for $\mathcal 
I(y_k, y_k, \epsilon^2)$, $\langle w_{jk} \rangle_{j 
\in J_k}$, with respect to the operator norm such that

\[ \#J_k \leq \left( \frac{1}{\epsilon} \right)^{C k^2 
\epsilon}. \]

\noindent Find a maximal sequence of $k \times k$ unitaries $u_{1k}, \ldots, u_{q_k k}$ such that for any $1 \leq i < j \leq q_k$,

\[ u_{ik} \Omega_k u_{ik}^*  \cap u_{jk} \Omega_k u_{jk}^* = \emptyset.  \] 

\noindent Consider the sets 
$\langle u_{jk} \Omega_k u_{jk}^* \rangle_{j=1}^{q_k}$.  By construction these 
sets are mutually disjoint and because there 
is a univeral lower bound on the volume of $\Omega_k$ and the sets are 
disjoint, $q_k$ is finite.  I will now find a lower bound on $q_k$.

Suppose $u \in U_k$.  Then by maximality of the $q_k$, for some $1 \leq i \leq q_k$ 

\[u \Omega_k u^* \cap u_{ik} \Omega_k u_{ik}^* \neq \emptyset.\] 

\noindent By definition of $\Omega_k$ and (6) this implies that for some 
$v, w \in B_k$, $u v\in u_{ik} w [\mathcal I (y_k, y_k, \epsilon^2)]^* 
\Rightarrow u \in u_{ik} w [ \mathcal I (y_k, y_k, \epsilon^2]^* v^* $.  
This being true for any $u \in U_k$,

\[U_k = \cup_{i=1}^{q_k} u_{ik} \mathcal N_{2 \epsilon, \infty} 
[\mathcal I (y_k, y_k, 4 \epsilon^2)]^*.\]   

\noindent For each $i$ there certainly exists a 
$3\epsilon$-cover for $u_{ik}\mathcal N_{2\epsilon, 
\infty}[\mathcal I(y_k,y_k, \epsilon^2)]$ with respect 
to the operator norm with no more than $\#J_k$ 
elements, and thus a $3 \epsilon$-cover for $U_k$ with 
no more that $q_k \cdot \#J_k$ elements.  Putting this 
together,

\[ \left( \frac{\alpha_0}{3\epsilon} \right)^{k^2} \leq 
q_k \cdot \#J_k \leq q_k \cdot \left( 
\frac{1}{\epsilon} \right)^{C k^2 \epsilon}. \]

\noindent Dividing on both sides,

\[ \left ( \frac{\alpha_0}{3} \right)^{k^2} \left(\frac{1}{\epsilon} 
\right)^{k^2(1- C\epsilon)}  \leq q_k.\]

\noindent As $\langle u_{jk} \Omega_k u_{jk}^* 
\rangle_{j=1}^{q_k}$ is a sequence of mutually disjoint 
subsets contained in the ball with operator norm $R$,

\[ q_k \cdot \text{vol} (\Omega_k) = \text{vol} \left 
(\cup_{j=1}^{q_k} u_{jk} \Omega_k u_{jk}^* \right) \leq 
C_{k,n} \cdot R^{nk^2} \]

\noindent and combining this with the preceding estimate,

\begin{eqnarray} \text{vol} (\Omega_k) \leq C_{k,n} 
R^{nk^2}(3/\alpha_0)^{k^2} \cdot 
\epsilon^{k^2 (1 - C \epsilon)}.\end{eqnarray}

\vspace{.1in}

\noindent Part 4: Finally, the estimate on $\text{vol}(\Omega_k)$ can be combined with Lemma 2.6 to get the desired upper bound on $\text{vol}(F_k)$.  
$\{uy_ku:u \in G_k\} \subset \Xi(y;m_1,k,\gamma_1)$ and 
by (5) and Lemma 2.6, for $k$ sufficiently large,

\begin{eqnarray} K_{\epsilon, \infty}(G_k) \leq \left ( 
\frac{1}{\epsilon} \right)^{\beta k^2} \cdot 
e^{k^2(\alpha_1 - \chi(y) + \log 5)}. \end{eqnarray}

\noindent The estimate above and (7) are enough to 
dominate $\text{vol}(F_k)$.  Suppose $\langle u_{lk} 
\rangle_{l \in \Lambda_k}$ is an $\epsilon$-net for 
$G_k$ with respect to the operator norm such that 
$\#\Lambda_k = K_{\epsilon, \infty}(G_k)$.  If $u \in 
G_k$ and $\xi \in L(m_1,k,\gamma_1)$, then for some $l 
\in \Lambda_k$, $\|u_{lk} - u\| < \epsilon$.  $u \xi 
u^* = u_{lk} (u_{lk}^* u \xi u^* u_{lk}) u_{lk}^*$ and 
$\|u_{lk}^* u
- I \| = \|u_{lk} - u\| < \epsilon$ so that $u_{lk}^* u 
\xi u^* u_{lk} \in \Omega_k$.  Thus, $u\xi u^* \in 
\cup_{l \in \Lambda} u_{lk} \Omega_k u_{lk}^* 
\Rightarrow \text{vol} (\cup_{u \in G_k} u 
L(m_1,k,\gamma_1) u^*) \leq K_{\epsilon, \infty}(G_k) 
\cdot \text{vol}(\Omega_k)$.  This with Part 2, (7), and (8) 
give:

\begin{eqnarray*} \text{vol}(\Xi_R(X;m_3,k,\gamma_3)) & 
\leq & \text{vol}(F_k) \\ & = & \text{vol}( \cup_{u \in 
G_k} u L(m_1,k,\gamma_1) u^* )\\ & \leq & K_{\epsilon, 
\infty}(G_k) \cdot \text{vol}(\Omega_k) \\ & \leq & 
\epsilon^{k^2(1 - L \epsilon - \beta)} \cdot C_{k,n} 
\cdot \left [
e^{k^2(\alpha_1 - \chi(y) + \log 5)} R^{nk^2} 
(3/\alpha_0)^{k^2} \right]    \\ \end{eqnarray*}

\noindent Sticking this estimate into the definition of 
$\chi_R(\Xi(X))$ yields

\[ \chi_R(\Xi(X)) \leq (1 - C \epsilon - \beta) \log 
\epsilon + \alpha_1 - \chi(y) + \log 5 + \log( 3 R^n 
\alpha_0^{-1}).  \]

\noindent This is true for any $0 < \epsilon < 
\epsilon_0$.  $C, R, \alpha_0, \alpha_1$ and $y$ 
are independent of $\epsilon$, 
and $\beta \in (0,1)$, so from all this it follows that 
the dominating term converges to $-\infty$ as $\epsilon 
\rightarrow 0$, whence $\chi_R(\Xi(X)) = -\infty$.  
$\chi(X \cup \{z\}) = - \infty$ as desired.

\end{proof}

\begin{corollary} If $y=y^*$ is in the von Neumann algebra generated by $X$, then $\chi(X \cup \{y\}) = -\infty$.
\end{corollary}

Taking the contrapositive of the main theorem produces:

\begin{corollary} If $\chi(X \cup \{z\}) > - \infty$, then for any element $x$ in the von Neumann algebra generated by $X$, $\delta_0(x, z) = \delta_0(x) + \delta_0(z)$.  In particular, if $x$ and $z$ are both diffuse, then $\delta_0(x, z) =2$.
\end{corollary}

\section{Non-microstates free entropy estimates}

 This last section deals with the claim concerning $\chi^*$, namely that if $y$ lies in the algebra generated by $X$, then $\chi^*(X \cup \{y\})= -\infty$.   More specifically, I'll show that if $y$ is in the algebra generated by $X$, then $\delta^*(X \cup \{y\}) < \delta^*(X) + 1$ and this will in turn imply that $\chi^*(X \cup \{y\}) = -\infty$.  I'll adhere to the nonmicrostate notation set forth in \cite{dvv:entropy4}.  First for a few preliminaries which are trivial modifications of the arguments in \cite{dvv:entropy4}.

\begin{lemma} If for each $1 \leq i \leq n$, $\beta_i \geq \alpha_i \geq 0$, then

\[ \chi^*(x_1 + \beta_1 s_1, \ldots, x_n + \beta_n s_n:B) \geq \chi^*(x_1 + \alpha_1 s_1, \ldots, x_n + \alpha_n s_n:B).
\]
\end{lemma}

\begin{proof} Suppose for each $1 \leq i \leq n$, $b_i \geq a_i >0$ are real numbers.  Find a semicircular family $\{s_1^{\prime}, \ldots, s_n^{\prime}\}$ free from $M$.  Using the free additive convolution formula for semicircular elements and the free Stam Inequality in \cite{dvv:entropy4}

\begin{eqnarray*} \Phi^*(x_1 + b_1 s_1, \ldots, x_n + b_n s_n) & = & \Phi^*(x_1 + a_1s_1 + \sqrt{b_1^2 - a_1^2} s_1^{\prime}, \ldots, x_n + a_n s_n + \sqrt{b_n^2 - a^2}s_n^{\prime}) \\ & \leq & \Phi^*(x_1 + a_1 s_1, \ldots, x_n + a_n s_n).
\end{eqnarray*}

\noindent So for $t >0$,

\[ -\Phi^*(x_1 + \beta_1 s_1 + \sqrt{t}s_1^{\prime}, \ldots, x_n + \beta_n s_1 + \sqrt{t}s_n^{\prime}) \geq -\Phi^*(x_1 + \alpha_1s_1 + \sqrt{t}s_1^{\prime}, \ldots, x_n + \alpha_n s_n + \sqrt{t}s_n^{\prime})
\]

\noindent The definition of $\chi^*$ now implies the desired result.
\end{proof}

\begin{corollary} If $\chi^*(X) > -\infty$, then $\delta^*(X) = n$.
\end{corollary}

\begin{proof}By definition $\delta^*(X) = n + \limsup_{\epsilon \rightarrow 0} \frac{\chi^*(X+ \epsilon S)}{|\log \epsilon|}$ and by the preceding lemma and the upper bound on $\chi^*$ in \cite{dvv:entropy4}   

\[ n \geq  n + \limsup_{\epsilon \rightarrow 0} \frac{\chi^*(X + \epsilon S)}{|\log \epsilon|} \geq 
n + \limsup_{\epsilon \rightarrow 0} \frac{\chi^*(X)}{|\log \epsilon|} =  n.\]
\end{proof}

Recall that \cite{cs:1} (see also a similar quantity in \cite{aagard:dt}) introduced a notion of dimension for $X$ given by the expression

\[ n - \liminf_{\epsilon \rightarrow 0} \epsilon^2 \Phi^*(X + \epsilon S).\]

\noindent I will denote this quantity by $\delta_1(X)$ as opposed to the original notation ($\delta^{\star}$) which uses a five-pointed star; my preference for this notation is simply that for the myopic it can be easily confused with $\delta^*$ which uses an asterik.  It is a consequence of \cite{cs:1} that $\delta^*(X) \leq \delta_1(X)$. 

\begin{lemma} If $y$ is a selfadjoint element contained in the algebra generated by $X$, then $\delta^*(X \cup \{y\}) \leq \delta_1(X \cup \{y\}) <  \delta_1(X) +1 \leq n+1$ and in particular $\chi^*(X \cup \{y\}) = - \infty$.
\end{lemma}

\begin{proof} Find a polynomial $f$ in $n$ noncommuting variables such that $f(X)=y$.  There exist two polynomials $g$ and $h$ such that

\[ f(X+ Z) = f(X) + g(X, Z) + h(X, Z)
\]

\noindent where $g$ is a polynomial in $2n$-noncommutative variables which is of homogeneous degree $1$ in the last $n$-variables, and $h$ is a polynomial in $2n$-noncommutative variables whose monomials terms all have degree greater than or equal to $2$ in its last $n$-variables.  In order to bound $\delta_1(X \cup \{y\})$ from above, it suffices to bound $\Phi^*(X + \epsilon S, y+ \epsilon s_{n+1})$ from below.

Fix $\epsilon >0$.   

\begin{eqnarray} \Phi^*(X + \epsilon S, y + \epsilon s_{n+1}) & \geq & \Phi^*(X + \epsilon S) + |J(y+ \epsilon s_{n+1}: X+ \epsilon S)|_2^2 \end{eqnarray}

\noindent By Proposition 3.9 \cite{dvv:entropy4}, 

\begin{eqnarray} J(y+ \epsilon s_{n+1}: X+\epsilon S) = \epsilon^{-1} \cdot E_{W^*(X+\epsilon S \cup \{y + \epsilon s_{n+1}\})}(s_{n+1}) \end{eqnarray}

\noindent Now $f(x_1 + \epsilon s_1, \ldots, x_n + \epsilon s_n) - (y+ \epsilon s_{n+1}) \in W^*(X + \epsilon S \cup \{ y+ \epsilon s_{n+1}\})$ and expanding the $f$ term

\begin{eqnarray*} f(x_1 + \epsilon s_1, \ldots, x_n + \epsilon s_n) - (y+ \epsilon s_{n+1}) & =& f(x_1, \ldots, x_n) + g(X, \epsilon S) + h(X, \epsilon S) - y - \epsilon s_{n+1} \\ & = & g(X, \epsilon S) + h(X, \epsilon S) - \epsilon s_{n+1}.\\
\end{eqnarray*}

\noindent Notice that $|g(X, \epsilon S)|_2 = \epsilon |g(X, S)|_2 \leq \epsilon L_1$ and $|h(X, \epsilon S)|_2 \leq \epsilon^2 L_2$ where $L_1$ and $L_2$ are constants dependent only on $g$ and $h$. Set 

\[\xi = \epsilon^{-1} (\epsilon s_{n+1} - g(X,\epsilon S)- h(X, \epsilon S)) \in W^*(X+ \epsilon S \cup \{ y+ \epsilon s_{n+1}\})\]

\noindent and observe that

\[ \|\xi\|_2 \leq 1 + L_1 + \epsilon L_2.
\]

\noindent Thus,

\begin{eqnarray*} \|E_{W^*(X+ \epsilon S \cup \{y + \epsilon s_{n+1}\})}(s_{n+1})\|_2 & \geq & \|\xi\|_2^{-1} \cdot  \varphi(s_{n+1} \xi) \\ & = & \|\xi\|_2^{-1} \\ & \geq & (1 + L_1 + \epsilon L_2)^{-1} \\
\end{eqnarray*}

\noindent and plugging this into (9) and using (10) along the way 

\begin{eqnarray*} \delta_1(X \cup \{y\}) & = & n+1 - \liminf_{\epsilon \rightarrow 0} \epsilon^2 \Phi^*(X + \epsilon S, y+ \epsilon s_{n+1}) \\ & \leq & n+1 - \liminf_{\epsilon \rightarrow 0}  \epsilon^2 \left [ \Phi^*(X+ \epsilon S) + \|J(y+ \epsilon s_{n+1}:X+ \epsilon S)\|_2^2 \right ] \\ & \leq & n+1 - \liminf_{\epsilon \rightarrow 0} \left [ \epsilon^2 \Phi^*(X+ \epsilon S) + \|E_{W*(X+ \epsilon S \cup \{y + \epsilon s_{n+1})}(s_{n+1})\|_2^2 \right ] \\ & \leq & n+1 - \liminf_{\epsilon \rightarrow 0} \left [ \epsilon^2 \Phi^*(X + \epsilon S) \right ] - (1+ L_1)^{-2} \\ & = & \delta_1(X) + 1 - (1+ L_1)^{-2}. \end{eqnarray*} 

\noindent The rest follows from Corollary 2.2.

\end{proof}

\noindent{\it Acknowledgments.} I thank Narutaka Ozawa and Terence Tao for several useful conversations.  I also thank Peter Pribik for the editorial comments.  

\bibliographystyle{amsplain}

\begin{thebibliography}{10}

\bibitem{aagard:dt} Lars Aagaard, \emph{The non-microstates free entropy 
dimension of any DT-operator is $2$}, J. Funct. Anal. 218 (2004), no. 1, 
176-205

\bibitem{sh} Serban Belinschii, Hari Bercovici, \emph{A property of free 
entropy}, Pacific J. Math. 211 (2003), no.1, 35-40. 

\bibitem{cs:1} Alain Connes, Dimitri Shlyakhtenko, \emph{$L^2$-homology 
for von Neumann algebras}, J. Reine Angew. Math. 586 (2005), 125-168.

\bibitem{ge:prime} Liming Ge, \emph{Applications of free entropy to finite 
von Neumann algebras, II}, Ann. Math. 147 (1998) 143-157.

\bibitem{gs} Liming Ge, Junhao Shen, \emph{On free entropy dimension of 
finite von Neumann algebras}, Geom. Funct. Anal. 12, (2002), no.3, 
546-566.

\bibitem{jung:hf} K. Jung, \emph{The free entropy dimension of hyperfinite 
von Neumann algebras}, Trans. of the Amer. Math. Soc. (2003), no.12, 
5053-5089

\bibitem{jung:hid} K. Jung, \emph{A hyperfinite inequality for 
$\delta_0$}, Proc. Amer. Math. Soc. 134 (2006), no.7, 2099-2108.

\bibitem{jung:s1b} K. Jung, \emph{Strongly 1-bounded von Neumann 
algebras}, to appear in Geom. Funct. Anal.

\bibitem{mattila} Pertti Mattila, \emph{Geometry of sets and measures in 
euclidean spaces} Cambridge Studies in Advanced Mathematics, v.44, 1995.

\bibitem{mehta} M. L. Mehta, \emph{Random Matrices, 3rd ed.} New York: Academic Press, 2004.

\bibitem{oneil:graph} Toby C. Oneil, \emph{Graphs of continuous functions 
from $R$ to $R$ are not purely unrectifiable}, Real Analysis Exchange, \textbf{26}, pp. 445-447.

\bibitem{raymond} J. Saint-Raymond, \emph{Le volume des ideaux 
d'operateurs classiques}, Studia Math. 80 (1984), 63-75.

\bibitem{schultz} Hanne Schultz, \emph{Semicircularity, Gaussianity, and 
Monotonicity of Entropy}, to appear in Journal of Operator Theory.

\bibitem{szarek} Stanislaw Szarek, \emph{Metric entropy of homogeneous 
spaces, Quantum Probability}, (Gdensk, 1997), Banach Center Publications, 
vol.43, Polish Academy of Science, Warsaw, 1998, pp. 395-410.

\bibitem{dvv:entropypower} Stanislaw Szarek, Dan Voiculescu, \emph{Volumes 
of restricted Minkowski sums and the free analogue of the entropy power 
inequality}, Comm. Math. Phys. 178 (1996), no.3, 563-570.

\bibitem{dvv:entropy2}
D.-V. Voiculescu, \emph{The analogues of entropy and of {Fisher's} information
  measure in free probability theory {II}}, Invent. Math. \textbf{118} 
  (1994), 411--440.

\bibitem{dvv:entropy3}
\bysame, \emph{The analogues of entropy and of {Fisher}'s information measure
  in free probability theory, {III}}, Geom. Funct. Anal.
  \textbf{6} (1996), 172--199.
  
\bibitem{dvv:entropy4}
\bysame, \emph{The analogues of entropy of {Fisher}'s information measure in free probability theory, {IV}}, Invent. Math. \textbf{132} (1998), 189-227.  

\bibitem{dvv:entropysurvey}
\bysame, \emph{Free entropy}, Bull. London Math. Soc. \textbf{34} (2002),
  no.~3, 257--278.
  
\bibitem{vn} John von Neumann, \emph{Approximative properties of matrices of 
high finite order} Portugaliae Math. 3, (1942) 1--62.

\end{thebibliography}

\providecommand{\bysame}{\leavevmode\hbox to3em{\hrulefill}\thinspace}

\end{document}